\documentclass[12pt,reqno]{amsart}
\usepackage{geometry}
\numberwithin{equation}{section}
\usepackage{amsmath,amsthm}
\usepackage{amssymb}
\usepackage{multicol}
\usepackage{mathtools}
\mathtoolsset{showonlyrefs}
\usepackage{bbm}

\usepackage{float}
\usepackage{mathrsfs}
 \usepackage{tikz}
\usepackage{diagbox}
\usepackage{booktabs}
\usepackage{enumitem}
\usepackage[sans]{dsfont} 
\usepackage{bbm}
\usepackage{changebar}
\usepackage[colorlinks]{hyperref}
\usepackage{a4wide}
\usepackage{physics}  

\usepackage{pst-node} 
\usepackage{tikz-cd}
\usetikzlibrary{arrows.meta}
\usepackage{subcaption}

\usepackage{tabstackengine}
\setstackEOL{;}
\setstackTAB{,}
\setstacktabbedgap{1ex}
\setstackgap{L}{1.0\normalbaselineskip}

\usepackage{tikz}
\usetikzlibrary{matrix}
\usetikzlibrary{backgrounds}
\usetikzlibrary{patterns,fadings}
\usetikzlibrary{arrows,decorations.pathmorphing}
\usetikzlibrary{calc}

\usepackage{pgfplots}\pgfplotsset{compat=1.18}

\usepackage{faktor}

\usepackage{tikz}
\usetikzlibrary{arrows,shapes,automata,petri,positioning,calc}
\tikzset{
    place/.style={
        circle,
        thick,
        draw=black,
        fill=gray!50,
        minimum size=20mm,
    },
        state/.style={
        circle,
        thick,
        draw=blue!75,
        fill=blue!20,
        minimum size=20mm,
    },
}

\tikzset{
    cross/.pic = {
    \draw[rotate = 45] (-0.2,0) -- (0.2,0);
    \draw[rotate = 45] (0,-0.2) -- (0, 0.2);
    }
}

\usepackage{ulem}
\usepackage{setspace}

\theoremstyle{plain}
\newtheorem{Th}{Theorem}[section]
\newtheorem{Lemma}[Th]{Lemma}

\newtheorem{Prop}[Th]{Proposition}

\theoremstyle{definition}
\newtheorem{Def}[Th]{Definition}

\newtheorem{?}[Th]{Problem}

\usepackage{soul} 

\newcommand{\rmd}{\mathrm{d}}   
\newcommand{\inner}[1]{{\langle #1\rangle}}
\newcommand{\nnorm}[1]{{\left\vert\kern-0.25ex\left\vert\kern-0.25ex\left\vert #1 
		\right\vert\kern-0.25ex\right\vert\kern-0.25ex\right\vert}} 

\usepackage{comment} 
\excludecomment{toexclude} 

\usepackage{stmaryrd} 

\usepackage{xfrac} 

\title[A gradient model for the Bernstein polynomial basis]{A gradient model for the Bernstein polynomial basis}

\author{Gabriel S. Nahum}
\address{Gabriel S. Nahum\\
inria Lyon, DRACULA team\\
56 Bd. Niels Bohr, 69199 Villeurbanne, France}
\email{{\tt gabriel.nahum@protonmail.com}}

\usepackage{multicol}
\begin{document}
\begin{small}
\textit{*The present version substantially simplifies the presentation and omits the proof of the hydrodynamic limit, which is now provided in a companion paper.}
\end{small}

\begin{abstract}
We introduce a symmetric, gradient exclusion process within the class of non-cooperative kinetically constrained lattice gases, modelling a non-linear diffusivity in which the exchange of occupation values between two neighbouring sites depends on the local density in specific boxes surrounding the pair. The existence of such a model satisfying the gradient property is the main novelty of this work, filling a gap in the literature regarding the types of diffusivities attainable within this class of models. The resulting dynamics exhibits similarities with the Bernstein polynomial basis and generalises the Porous Media Model. We also introduce an auxiliary collection of processes, which extend the Porous Media Model in a different direction and are related to the former process via an inversion formula.
\end{abstract}

\maketitle
\begingroup
\renewcommand{\thefootnote}{}

\footnotetext{\textbf{Acknowledgements:} The author warmly thanks Marielle Simon and Patrícia Gonçalves for their insightful suggestions and support throughout this work, which greatly enriched this study.}

\addtocounter{footnote}{-1}
\endgroup




\tableofcontents

\section{Introduction}

There is an ongoing effort to rigorously understand macroscopic behaviour and thermodynamic properties in non-equilibrium statistical mechanics through the use of Markovian interacting particle systems to model the evolution of physical systems on a microscopic scale. In the context of symmetric exclusion processes, and at the core of this work, the Porous Media Model \cite{GLT} (PMM) is a continuous-time Markov process with a rich dynamics, describing degenerate diffusion. Initially, particles are randomly distributed on a discrete unidimensional lattice under the \textit{exclusion rule}, meaning that at most one particle may occupy a site. Each pair of neighbouring sites is independently associated with a Poisson process, dictating the times of potential jumps, with specific rates depending on the local configuration in the following way: fixing some natural number $n$, a pair of neighbouring sites exchanges their occupation values if there are $n$ aligned particles around the pair, with a rate given by the proportion of such groups of aligned particles. Throughout this work, for simplicity we fix Let $N\in\mathbb{N}_+$ be fixed and consider the discrete lattice given by $\mathbb{T}_N:=\mathbb{T}/N\mathbb{Z}$, the one-dimensional discrete torus. The dynamics is illustrated in Figure \ref{fig:PMM}. A motivation for introducing this process is to model the degeneracy of diffusion at trivial density $\rho\to 0$, within the class of \textit{non-cooperative kinetically constrained lattice gases}. The "non-cooperative" aspect refers to the existence of \textit{mobile clusters}—clusters of particles that, through individual jumps, can perform excursions in the lattice. Notably, it has been formalized in various settings that the local density of particles is associated with the Porous Media Equation \cite{GLT,BDGN,CG24,s2f,CGN24} with diffusion coefficient $D(\rho)=\rho^n$.

\begin{figure}[H]
\centering

 
\tikzset{
pattern size/.store in=\mcSize, 
pattern size = 5pt,
pattern thickness/.store in=\mcThickness, 
pattern thickness = 0.3pt,
pattern radius/.store in=\mcRadius, 
pattern radius = 1pt}
\makeatletter
\pgfutil@ifundefined{pgf@pattern@name@_bekc5s796}{
\pgfdeclarepatternformonly[\mcThickness,\mcSize]{_bekc5s796}
{\pgfqpoint{0pt}{-\mcThickness}}
{\pgfpoint{\mcSize}{\mcSize}}
{\pgfpoint{\mcSize}{\mcSize}}
{
\pgfsetcolor{\tikz@pattern@color}
\pgfsetlinewidth{\mcThickness}
\pgfpathmoveto{\pgfqpoint{0pt}{\mcSize}}
\pgfpathlineto{\pgfpoint{\mcSize+\mcThickness}{-\mcThickness}}
\pgfusepath{stroke}
}}
\makeatother

 
\tikzset{
pattern size/.store in=\mcSize, 
pattern size = 5pt,
pattern thickness/.store in=\mcThickness, 
pattern thickness = 0.3pt,
pattern radius/.store in=\mcRadius, 
pattern radius = 1pt}
\makeatletter
\pgfutil@ifundefined{pgf@pattern@name@_8phuy4aw0}{
\pgfdeclarepatternformonly[\mcThickness,\mcSize]{_8phuy4aw0}
{\pgfqpoint{0pt}{-\mcThickness}}
{\pgfpoint{\mcSize}{\mcSize}}
{\pgfpoint{\mcSize}{\mcSize}}
{
\pgfsetcolor{\tikz@pattern@color}
\pgfsetlinewidth{\mcThickness}
\pgfpathmoveto{\pgfqpoint{0pt}{\mcSize}}
\pgfpathlineto{\pgfpoint{\mcSize+\mcThickness}{-\mcThickness}}
\pgfusepath{stroke}
}}
\makeatother
\tikzset{every picture/.style={line width=0.75pt}} 

\begin{tikzpicture}[x=0.75pt,y=0.75pt,yscale=-1,xscale=1]

\draw    (121.5,219.58) -- (494.5,220.58) (155.51,215.67) -- (155.49,223.67)(189.51,215.76) -- (189.49,223.76)(223.51,215.85) -- (223.49,223.85)(257.51,215.94) -- (257.49,223.94)(291.51,216.04) -- (291.49,224.04)(325.51,216.13) -- (325.49,224.13)(359.51,216.22) -- (359.49,224.22)(393.51,216.31) -- (393.49,224.31)(427.51,216.4) -- (427.49,224.4)(461.51,216.49) -- (461.49,224.49) ;
\draw  [fill={rgb, 255:red, 155; green, 155; blue, 155 }  ,fill opacity=1 ] (277.14,199.35) .. controls (277.14,190.24) and (284.73,182.85) .. (294.09,182.85) .. controls (303.45,182.85) and (311.05,190.24) .. (311.05,199.35) .. controls (311.05,208.46) and (303.45,215.85) .. (294.09,215.85) .. controls (284.73,215.85) and (277.14,208.46) .. (277.14,199.35) -- cycle ;
\draw   (311.05,199.35) .. controls (311.05,190.24) and (318.64,182.85) .. (328,182.85) .. controls (337.36,182.85) and (344.95,190.24) .. (344.95,199.35) .. controls (344.95,208.46) and (337.36,215.85) .. (328,215.85) .. controls (318.64,215.85) and (311.05,208.46) .. (311.05,199.35) -- cycle ;
\draw  [fill={rgb, 255:red, 155; green, 155; blue, 155 }  ,fill opacity=1 ] (344.95,199.35) .. controls (344.95,190.24) and (352.55,182.85) .. (361.91,182.85) .. controls (371.27,182.85) and (378.86,190.24) .. (378.86,199.35) .. controls (378.86,208.46) and (371.27,215.85) .. (361.91,215.85) .. controls (352.55,215.85) and (344.95,208.46) .. (344.95,199.35) -- cycle ;
\draw  [fill={rgb, 255:red, 155; green, 155; blue, 155 }  ,fill opacity=1 ] (378.86,199.35) .. controls (378.86,190.24) and (386.45,182.85) .. (395.82,182.85) .. controls (405.18,182.85) and (412.77,190.24) .. (412.77,199.35) .. controls (412.77,208.46) and (405.18,215.85) .. (395.82,215.85) .. controls (386.45,215.85) and (378.86,208.46) .. (378.86,199.35) -- cycle ;
\draw  [fill={rgb, 255:red, 155; green, 155; blue, 155 }  ,fill opacity=1 ] (412.77,199.35) .. controls (412.77,190.24) and (420.36,182.85) .. (429.73,182.85) .. controls (439.09,182.85) and (446.68,190.24) .. (446.68,199.35) .. controls (446.68,208.46) and (439.09,215.85) .. (429.73,215.85) .. controls (420.36,215.85) and (412.77,208.46) .. (412.77,199.35) -- cycle ;
\draw  [fill={rgb, 255:red, 155; green, 155; blue, 155 }  ,fill opacity=1 ] (243.23,199.35) .. controls (243.23,190.24) and (250.82,182.85) .. (260.18,182.85) .. controls (269.55,182.85) and (277.14,190.24) .. (277.14,199.35) .. controls (277.14,208.46) and (269.55,215.85) .. (260.18,215.85) .. controls (250.82,215.85) and (243.23,208.46) .. (243.23,199.35) -- cycle ;
\draw  [fill={rgb, 255:red, 155; green, 155; blue, 155 }  ,fill opacity=1 ] (209.32,199.35) .. controls (209.32,190.24) and (216.91,182.85) .. (226.27,182.85) .. controls (235.64,182.85) and (243.23,190.24) .. (243.23,199.35) .. controls (243.23,208.46) and (235.64,215.85) .. (226.27,215.85) .. controls (216.91,215.85) and (209.32,208.46) .. (209.32,199.35) -- cycle ;
\draw   (175.41,199.35) .. controls (175.41,190.24) and (183,182.85) .. (192.36,182.85) .. controls (201.73,182.85) and (209.32,190.24) .. (209.32,199.35) .. controls (209.32,208.46) and (201.73,215.85) .. (192.36,215.85) .. controls (183,215.85) and (175.41,208.46) .. (175.41,199.35) -- cycle ;
\draw   (141.5,199.35) .. controls (141.5,190.24) and (149.09,182.85) .. (158.45,182.85) .. controls (167.82,182.85) and (175.41,190.24) .. (175.41,199.35) .. controls (175.41,208.46) and (167.82,215.85) .. (158.45,215.85) .. controls (149.09,215.85) and (141.5,208.46) .. (141.5,199.35) -- cycle ;
\draw   (446.68,199.35) .. controls (446.68,190.24) and (454.27,182.85) .. (463.64,182.85) .. controls (473,182.85) and (480.59,190.24) .. (480.59,199.35) .. controls (480.59,208.46) and (473,215.85) .. (463.64,215.85) .. controls (454.27,215.85) and (446.68,208.46) .. (446.68,199.35) -- cycle ;
\draw  [color={rgb, 255:red, 0; green, 0; blue, 0 }  ,draw opacity=1 ][pattern=_bekc5s796,pattern size=6pt,pattern thickness=0.75pt,pattern radius=0pt, pattern color={rgb, 255:red, 0; green, 0; blue, 0}][line width=1.5]  (243.23,116.85) -- (446.68,116.85) -- (446.68,149.85) -- (243.23,149.85) -- cycle ;
\draw  [color={rgb, 255:red, 0; green, 0; blue, 0 }  ,draw opacity=1 ][pattern=_8phuy4aw0,pattern size=6pt,pattern thickness=0.75pt,pattern radius=0pt, pattern color={rgb, 255:red, 0; green, 0; blue, 0}][line width=1.5]  (209.32,149.85) -- (412.77,149.85) -- (412.77,182.85) -- (209.32,182.85) -- cycle ;
\draw  [color={rgb, 255:red, 208; green, 2; blue, 27 }  ,draw opacity=1 ][line width=1.5]  (277.14,182.85) -- (344.95,182.85) -- (344.95,215.85) -- (277.14,215.85) -- cycle ;

\draw (487,191.09) node [anchor=north west][inner sep=0.75pt]    {$\cdots $};
\draw (108.73,191.09) node [anchor=north west][inner sep=0.75pt]    {$\cdots $};
\draw (310.17,229.3) node [anchor=north west][inner sep=0.75pt]  [font=\tiny] [align=left] {$\displaystyle x+1$};
\draw (285,231) node [anchor=north west][inner sep=0.75pt]  [font=\tiny] [align=left] {$\displaystyle x$};

\end{tikzpicture}    
\caption{PMM rate ($n=4$) for an exchange in the node $\{x,x+1\}$. The patterned rectangles represent the boxes where there are $n$ aligned particles around $\{x,x+1\}$. The total rate is $2/5$.}\label{fig:PMM}
 \end{figure}

A convenient characteristic of the PMM, and one relevant to the present work, is that it enjoys the gradient property. This means precisely that the algebraic microscopic current can be expressed as the (discrete) gradient of some function \cite[Part II Subsection 2.4]{spohn:book}. In this case, the diffusion coefficient of the hydrodynamic equation has an explicit expression, and a scaling limit analysis is much less involved than in the general setting \cite{phd:clement}. In the scaling-limit for a large variety of models, it is known that the bulk diffusion coefficient can be formally expressed by a variational problem \cite[Part II Subsection 2.2]{spohn:book}, and it is known that the gradient property is satisfied if \cite[Part II Subsection 2.4]{spohn:book} and only if \cite{makiko} the dynamical part of the aforementioned variational formula is zero, providing this property a physical meaning.

The motivation for the present work is the ongoing effort to understand what types of differential equations can be associated with gradient models, in particular the works  \cite{s2f,CGN24}. To be more precise, through the superposition of PMMs, one can derive diffusion coefficients $ D(\rho) $ that can be expressed in a polynomial basis, $ D(\rho) = a_0 + a_1\rho + a_2\rho^2 + \cdots $. Naturally, for the underlying model to be well-defined, the associated constraints must be non-negative, which imposes conditions on the coefficients $ (a_i)_{i\geq 0} $. We identified that it is not possible to construct a coefficient $ D(\rho) $ with two roots of multiplicity higher than two using the PMMs as a basis. In particular, the collection of PMMs cannot be identified with a monomial basis for the set of non-negative polynomials. The key example of an unattainable diffusivity is an element of the Bernstein polynomial basis, $ \text{B}_{n,L}(\rho) = \binom{L}{n} \rho^n (1-\rho)^{L-n} $, with $ 2 \leq n \leq L-2 $. The underlying reason is related to the monotonicity of the sequence of constraints of the PMM, as we explain later on. This leads to the following questions:  
\begin{enumerate}  
    \item Is there a gradient dynamics associated with these types of diffusion coefficients?  
    \item If it exists, what is the dynamics associated with the underlying monomial basis $ \{\rho^n\}_{n\geq 0} $, and how does it relate to the PMM?  
\end{enumerate}  
To fill this gap in the literature, we focus on understanding how to generalize the PMM dynamics in such a way that at the diffusivity level, $\rho^L\mapsto\rho^n(1-\rho)^{L-n}$, while also maintaining the gradient property.

The answer to the first question is the following dynamics: given fixed natural numbers $ n \leq L $, a pair composed of an empty site and an occupied site exchanges occupation values at a rate given by the proportion of boxes of length $ L+2 $ with exactly $ n+1 $ particles that contain the pair of neighbouring sites. Consequently, boxes of length $ L+2 $ with exactly $ n+1 $ particles are \textit{mobile clusters}. The dynamics is exemplified in Figure \ref{fig:PMMnk} and will be referred to as the Bernstein model.

\begin{figure}[H]
\centering

 
\tikzset{
pattern size/.store in=\mcSize, 
pattern size = 5pt,
pattern thickness/.store in=\mcThickness, 
pattern thickness = 0.3pt,
pattern radius/.store in=\mcRadius, 
pattern radius = 1pt}
\makeatletter
\pgfutil@ifundefined{pgf@pattern@name@_dnwmrzisr}{
\pgfdeclarepatternformonly[\mcThickness,\mcSize]{_dnwmrzisr}
{\pgfqpoint{0pt}{0pt}}
{\pgfpoint{\mcSize+\mcThickness}{\mcSize+\mcThickness}}
{\pgfpoint{\mcSize}{\mcSize}}
{
\pgfsetcolor{\tikz@pattern@color}
\pgfsetlinewidth{\mcThickness}
\pgfpathmoveto{\pgfqpoint{0pt}{0pt}}
\pgfpathlineto{\pgfpoint{\mcSize+\mcThickness}{\mcSize+\mcThickness}}
\pgfusepath{stroke}
}}
\makeatother
\tikzset{every picture/.style={line width=0.75pt}} 

\begin{tikzpicture}[x=0.75pt,y=0.75pt,yscale=-1,xscale=1]

\draw    (141.5,239.58) -- (514.5,240.58) (175.51,235.67) -- (175.49,243.67)(209.51,235.76) -- (209.49,243.76)(243.51,235.85) -- (243.49,243.85)(277.51,235.94) -- (277.49,243.94)(311.51,236.04) -- (311.49,244.04)(345.51,236.13) -- (345.49,244.13)(379.51,236.22) -- (379.49,244.22)(413.51,236.31) -- (413.49,244.31)(447.51,236.4) -- (447.49,244.4)(481.51,236.49) -- (481.49,244.49) ;
\draw  [fill={rgb, 255:red, 155; green, 155; blue, 155 }  ,fill opacity=1 ] (297.14,219.35) .. controls (297.14,210.24) and (304.73,202.85) .. (314.09,202.85) .. controls (323.45,202.85) and (331.05,210.24) .. (331.05,219.35) .. controls (331.05,228.46) and (323.45,235.85) .. (314.09,235.85) .. controls (304.73,235.85) and (297.14,228.46) .. (297.14,219.35) -- cycle ;
\draw   (331.05,219.35) .. controls (331.05,210.24) and (338.64,202.85) .. (348,202.85) .. controls (357.36,202.85) and (364.95,210.24) .. (364.95,219.35) .. controls (364.95,228.46) and (357.36,235.85) .. (348,235.85) .. controls (338.64,235.85) and (331.05,228.46) .. (331.05,219.35) -- cycle ;
\draw  [fill={rgb, 255:red, 155; green, 155; blue, 155 }  ,fill opacity=1 ] (364.95,219.35) .. controls (364.95,210.24) and (372.55,202.85) .. (381.91,202.85) .. controls (391.27,202.85) and (398.86,210.24) .. (398.86,219.35) .. controls (398.86,228.46) and (391.27,235.85) .. (381.91,235.85) .. controls (372.55,235.85) and (364.95,228.46) .. (364.95,219.35) -- cycle ;
\draw  [fill={rgb, 255:red, 155; green, 155; blue, 155 }  ,fill opacity=1 ] (398.86,219.35) .. controls (398.86,210.24) and (406.45,202.85) .. (415.82,202.85) .. controls (425.18,202.85) and (432.77,210.24) .. (432.77,219.35) .. controls (432.77,228.46) and (425.18,235.85) .. (415.82,235.85) .. controls (406.45,235.85) and (398.86,228.46) .. (398.86,219.35) -- cycle ;
\draw  [fill={rgb, 255:red, 155; green, 155; blue, 155 }  ,fill opacity=1 ] (432.77,219.35) .. controls (432.77,210.24) and (440.36,202.85) .. (449.73,202.85) .. controls (459.09,202.85) and (466.68,210.24) .. (466.68,219.35) .. controls (466.68,228.46) and (459.09,235.85) .. (449.73,235.85) .. controls (440.36,235.85) and (432.77,228.46) .. (432.77,219.35) -- cycle ;
\draw  [fill={rgb, 255:red, 155; green, 155; blue, 155 }  ,fill opacity=1 ] (263.23,219.35) .. controls (263.23,210.24) and (270.82,202.85) .. (280.18,202.85) .. controls (289.55,202.85) and (297.14,210.24) .. (297.14,219.35) .. controls (297.14,228.46) and (289.55,235.85) .. (280.18,235.85) .. controls (270.82,235.85) and (263.23,228.46) .. (263.23,219.35) -- cycle ;
\draw  [fill={rgb, 255:red, 155; green, 155; blue, 155 }  ,fill opacity=1 ] (229.32,219.35) .. controls (229.32,210.24) and (236.91,202.85) .. (246.27,202.85) .. controls (255.64,202.85) and (263.23,210.24) .. (263.23,219.35) .. controls (263.23,228.46) and (255.64,235.85) .. (246.27,235.85) .. controls (236.91,235.85) and (229.32,228.46) .. (229.32,219.35) -- cycle ;
\draw   (195.41,219.35) .. controls (195.41,210.24) and (203,202.85) .. (212.36,202.85) .. controls (221.73,202.85) and (229.32,210.24) .. (229.32,219.35) .. controls (229.32,228.46) and (221.73,235.85) .. (212.36,235.85) .. controls (203,235.85) and (195.41,228.46) .. (195.41,219.35) -- cycle ;
\draw   (161.5,219.35) .. controls (161.5,210.24) and (169.09,202.85) .. (178.45,202.85) .. controls (187.82,202.85) and (195.41,210.24) .. (195.41,219.35) .. controls (195.41,228.46) and (187.82,235.85) .. (178.45,235.85) .. controls (169.09,235.85) and (161.5,228.46) .. (161.5,219.35) -- cycle ;
\draw   (466.68,219.35) .. controls (466.68,210.24) and (474.27,202.85) .. (483.64,202.85) .. controls (493,202.85) and (500.59,210.24) .. (500.59,219.35) .. controls (500.59,228.46) and (493,235.85) .. (483.64,235.85) .. controls (474.27,235.85) and (466.68,228.46) .. (466.68,219.35) -- cycle ;
\draw  [color={rgb, 255:red, 0; green, 0; blue, 0 }  ,draw opacity=1 ][pattern=_dnwmrzisr,pattern size=6pt,pattern thickness=0.75pt,pattern radius=0pt, pattern color={rgb, 255:red, 0; green, 0; blue, 0}][line width=1.5]  (161.5,169.85) -- (364.95,169.85) -- (364.95,202.85) -- (161.5,202.85) -- cycle ;
\draw  [color={rgb, 255:red, 208; green, 2; blue, 27 }  ,draw opacity=1 ][line width=1.5]  (297.14,202.85) -- (364.95,202.85) -- (364.95,235.85) -- (297.14,235.85) -- cycle ;

\draw (508.5,211.09) node [anchor=north west][inner sep=0.75pt]    {$\cdots $};
\draw (128.73,211.09) node [anchor=north west][inner sep=0.75pt]    {$\cdots $};
\draw (333,249) node [anchor=north west][inner sep=0.75pt]  [font=\tiny] [align=left] {$x+1$};
\draw (307,251.2) node [anchor=north west][inner sep=0.75pt]  [font=\tiny] [align=left] {$x$};

\end{tikzpicture}\caption{Bernstein model ($n=2,L=4$) rate for an exchange in the node $\{x,x+1\}$. The patterned rectangle represents the box containing $\{x,x+1\}$ with exactly $n+1$ particles. The total rate is $1/5$.}\label{fig:PMMnk}
\end{figure}

Remarkably, this simple dynamics enjoys the gradient property and generalises the combinatorial mechanism of the PMM, while coinciding with a PMM for $ n=L $. Moreover, it is associated with the diffusion coefficient $ D(\rho) = \text{B}_{n,L}(\rho) $. The process also enjoys mobile clusters that generalise those of the PMM, and algebraic properties analogous to the Bernstein polynomial basis, which also relate it to the Symmetric Simple Exclusion Process.  

Regarding question (2), the dynamics is as follows. Fixed $ 0\leq \ell\leq L $, a pair composed of two nearest-neighbour sites, one occupied and one vacant, exchanges its occupation values if it is contained in at least one box of length $ L+2 $ with at least $ \ell+1 $ particles and one vacant site. The precise rate of exchange is given by the uniform average over the possible boxes of length $ L+2 $, of the probability of choosing $ \ell $ occupied sites independently out of the existing $ m\geq \ell $ occupied sites in the respective $ L+2 $-length box. This dynamics also extends the PMM, coinciding with one for $ \ell=L $, and being associated with the diffusion coefficient $ D(\rho) = \rho^\ell $ for each $ 0\leq \ell\leq L $. We refer to this process as the \textit{reduced} Porous Media Model since it corresponds, fixed $L$, to a reduction of the diffusivity of the PMM, $ \alpha^L\to\alpha^\ell $, while maintaining the range of interaction, $ 2L+2 $. The dynamics is illustrated in Figure \ref{fig:redPMM} below.

\begin{figure}[H]
\centering

 
\tikzset{
pattern size/.store in=\mcSize, 
pattern size = 5pt,
pattern thickness/.store in=\mcThickness, 
pattern thickness = 0.3pt,
pattern radius/.store in=\mcRadius, 
pattern radius = 1pt}
\makeatletter
\pgfutil@ifundefined{pgf@pattern@name@_9biovjtxw}{
\pgfdeclarepatternformonly[\mcThickness,\mcSize]{_9biovjtxw}
{\pgfqpoint{-\mcThickness}{-\mcThickness}}
{\pgfpoint{\mcSize}{\mcSize}}
{\pgfpoint{\mcSize}{\mcSize}}
{
\pgfsetcolor{\tikz@pattern@color}
\pgfsetlinewidth{\mcThickness}
\pgfpathmoveto{\pgfpointorigin}
\pgfpathlineto{\pgfpoint{0}{\mcSize}}
\pgfusepath{stroke}
}}
\makeatother

 
\tikzset{
pattern size/.store in=\mcSize, 
pattern size = 5pt,
pattern thickness/.store in=\mcThickness, 
pattern thickness = 0.3pt,
pattern radius/.store in=\mcRadius, 
pattern radius = 1pt}
\makeatletter
\pgfutil@ifundefined{pgf@pattern@name@_9ahfwx93n}{
\pgfdeclarepatternformonly[\mcThickness,\mcSize]{_9ahfwx93n}
{\pgfqpoint{0pt}{0pt}}
{\pgfpoint{\mcSize+\mcThickness}{\mcSize+\mcThickness}}
{\pgfpoint{\mcSize}{\mcSize}}
{
\pgfsetcolor{\tikz@pattern@color}
\pgfsetlinewidth{\mcThickness}
\pgfpathmoveto{\pgfqpoint{0pt}{0pt}}
\pgfpathlineto{\pgfpoint{\mcSize+\mcThickness}{\mcSize+\mcThickness}}
\pgfusepath{stroke}
}}
\makeatother

 
\tikzset{
pattern size/.store in=\mcSize, 
pattern size = 5pt,
pattern thickness/.store in=\mcThickness, 
pattern thickness = 0.3pt,
pattern radius/.store in=\mcRadius, 
pattern radius = 1pt}
\makeatletter
\pgfutil@ifundefined{pgf@pattern@name@_tt1aga6dd}{
\pgfdeclarepatternformonly[\mcThickness,\mcSize]{_tt1aga6dd}
{\pgfqpoint{0pt}{-\mcThickness}}
{\pgfpoint{\mcSize}{\mcSize}}
{\pgfpoint{\mcSize}{\mcSize}}
{
\pgfsetcolor{\tikz@pattern@color}
\pgfsetlinewidth{\mcThickness}
\pgfpathmoveto{\pgfqpoint{0pt}{\mcSize}}
\pgfpathlineto{\pgfpoint{\mcSize+\mcThickness}{-\mcThickness}}
\pgfusepath{stroke}
}}
\makeatother
\tikzset{every picture/.style={line width=0.75pt}} 

\begin{tikzpicture}[x=0.75pt,y=0.75pt,yscale=-1,xscale=1]

\draw  [color={rgb, 255:red, 0; green, 0; blue, 0 }  ,draw opacity=1 ][pattern=_9biovjtxw,pattern size=4.5pt,pattern thickness=0.75pt,pattern radius=0pt, pattern color={rgb, 255:red, 0; green, 0; blue, 0}][line width=1.5]  (209.32,83.85) -- (412.77,83.85) -- (412.77,116.85) -- (209.32,116.85) -- cycle ;
\draw    (121.5,219.58) -- (494.5,220.58) (155.51,215.67) -- (155.49,223.67)(189.51,215.76) -- (189.49,223.76)(223.51,215.85) -- (223.49,223.85)(257.51,215.94) -- (257.49,223.94)(291.51,216.04) -- (291.49,224.04)(325.51,216.13) -- (325.49,224.13)(359.51,216.22) -- (359.49,224.22)(393.51,216.31) -- (393.49,224.31)(427.51,216.4) -- (427.49,224.4)(461.51,216.49) -- (461.49,224.49) ;
\draw  [fill={rgb, 255:red, 155; green, 155; blue, 155 }  ,fill opacity=1 ] (277.14,199.35) .. controls (277.14,190.24) and (284.73,182.85) .. (294.09,182.85) .. controls (303.45,182.85) and (311.05,190.24) .. (311.05,199.35) .. controls (311.05,208.46) and (303.45,215.85) .. (294.09,215.85) .. controls (284.73,215.85) and (277.14,208.46) .. (277.14,199.35) -- cycle ;
\draw   (311.05,199.35) .. controls (311.05,190.24) and (318.64,182.85) .. (328,182.85) .. controls (337.36,182.85) and (344.95,190.24) .. (344.95,199.35) .. controls (344.95,208.46) and (337.36,215.85) .. (328,215.85) .. controls (318.64,215.85) and (311.05,208.46) .. (311.05,199.35) -- cycle ;
\draw   (344.95,199.35) .. controls (344.95,190.24) and (352.55,182.85) .. (361.91,182.85) .. controls (371.27,182.85) and (378.86,190.24) .. (378.86,199.35) .. controls (378.86,208.46) and (371.27,215.85) .. (361.91,215.85) .. controls (352.55,215.85) and (344.95,208.46) .. (344.95,199.35) -- cycle ;
\draw   (378.86,199.35) .. controls (378.86,190.24) and (386.45,182.85) .. (395.82,182.85) .. controls (405.18,182.85) and (412.77,190.24) .. (412.77,199.35) .. controls (412.77,208.46) and (405.18,215.85) .. (395.82,215.85) .. controls (386.45,215.85) and (378.86,208.46) .. (378.86,199.35) -- cycle ;
\draw   (412.77,199.35) .. controls (412.77,190.24) and (420.36,182.85) .. (429.73,182.85) .. controls (439.09,182.85) and (446.68,190.24) .. (446.68,199.35) .. controls (446.68,208.46) and (439.09,215.85) .. (429.73,215.85) .. controls (420.36,215.85) and (412.77,208.46) .. (412.77,199.35) -- cycle ;
\draw  [fill={rgb, 255:red, 155; green, 155; blue, 155 }  ,fill opacity=1 ] (243.23,199.35) .. controls (243.23,190.24) and (250.82,182.85) .. (260.18,182.85) .. controls (269.55,182.85) and (277.14,190.24) .. (277.14,199.35) .. controls (277.14,208.46) and (269.55,215.85) .. (260.18,215.85) .. controls (250.82,215.85) and (243.23,208.46) .. (243.23,199.35) -- cycle ;
\draw  [fill={rgb, 255:red, 155; green, 155; blue, 155 }  ,fill opacity=1 ] (209.32,199.35) .. controls (209.32,190.24) and (216.91,182.85) .. (226.27,182.85) .. controls (235.64,182.85) and (243.23,190.24) .. (243.23,199.35) .. controls (243.23,208.46) and (235.64,215.85) .. (226.27,215.85) .. controls (216.91,215.85) and (209.32,208.46) .. (209.32,199.35) -- cycle ;
\draw   (175.41,199.35) .. controls (175.41,190.24) and (183,182.85) .. (192.36,182.85) .. controls (201.73,182.85) and (209.32,190.24) .. (209.32,199.35) .. controls (209.32,208.46) and (201.73,215.85) .. (192.36,215.85) .. controls (183,215.85) and (175.41,208.46) .. (175.41,199.35) -- cycle ;
\draw   (141.5,199.35) .. controls (141.5,190.24) and (149.09,182.85) .. (158.45,182.85) .. controls (167.82,182.85) and (175.41,190.24) .. (175.41,199.35) .. controls (175.41,208.46) and (167.82,215.85) .. (158.45,215.85) .. controls (149.09,215.85) and (141.5,208.46) .. (141.5,199.35) -- cycle ;
\draw   (446.68,199.35) .. controls (446.68,190.24) and (454.27,182.85) .. (463.64,182.85) .. controls (473,182.85) and (480.59,190.24) .. (480.59,199.35) .. controls (480.59,208.46) and (473,215.85) .. (463.64,215.85) .. controls (454.27,215.85) and (446.68,208.46) .. (446.68,199.35) -- cycle ;
\draw  [color={rgb, 255:red, 0; green, 0; blue, 0 }  ,draw opacity=1 ][pattern=_9ahfwx93n,pattern size=6pt,pattern thickness=0.75pt,pattern radius=0pt, pattern color={rgb, 255:red, 0; green, 0; blue, 0}][line width=1.5]  (141.5,149.85) -- (344.95,149.85) -- (344.95,182.85) -- (141.5,182.85) -- cycle ;
\draw  [color={rgb, 255:red, 208; green, 2; blue, 27 }  ,draw opacity=1 ][line width=1.5]  (277.14,182.85) -- (344.95,182.85) -- (344.95,215.85) -- (277.14,215.85) -- cycle ;
\draw  [color={rgb, 255:red, 0; green, 0; blue, 0 }  ,draw opacity=1 ][pattern=_tt1aga6dd,pattern size=6pt,pattern thickness=0.75pt,pattern radius=0pt, pattern color={rgb, 255:red, 0; green, 0; blue, 0}][line width=1.5]  (175.41,116.85) -- (378.86,116.85) -- (378.86,149.85) -- (175.41,149.85) -- cycle ;

\draw (486,191.09) node [anchor=north west][inner sep=0.75pt]    {$\cdots $};
\draw (108.73,191.09) node [anchor=north west][inner sep=0.75pt]    {$\cdots $};
\draw (311,229.3) node [anchor=north west][inner sep=0.75pt]  [font=\tiny] [align=left] {$x+1$};
\draw (286,231) node [anchor=north west][inner sep=0.75pt]  [font=\tiny] [align=left] {$x$};

\end{tikzpicture}
\caption{Reduced PMM ($\ell=2,L=4$) rates for a jump in $\{x,x+1\}$. The patterned rectangles represent the boxes containing $\{x,x+1\}$ with at least two particles. The total rate is $\tfrac{3}{5\binom{4}{2}}$.}\label{fig:redPMM}
\end{figure}

The two previously introduced models are connected through an inversion formula. Alternatively, introducing the dynamics from question (2), one can see the Bernstein model's constraints as a binomial transformation of a sequence of reduced PMMs' constraints.  

The content of this work is to prove the gradient property of the main models, connect them through an inversion formula, and show other relevant properties, thereby paving the way for a scaling-limit analysis in subsequent works. We note that although the main motivation lies in the context of statistical mechanics, our results and arguments are of a combinatorial nature.

\subsection{Future work}\label{subsec:future}

The models introduced here serve as a starting step for further developments in the study of gradient dynamics and their associated diffusivities.

A motivation for the introduction of the Bernstein model is the work in \cite{s2f}, where the authors extend the PMM into a family continuously parametrized by $m\in(0,2]$ and derive the equation $\partial_t\rho=\partial_u^2\rho^m$, thereby describing the transition from slow ($m>1$) to fast ($m<1$) diffusion. A tool missing for the extension of the dynamics to $m>2$ is precisely the collection of models introduced here. This is performed in a companion paper, requiring novel arguments at the level of the so-called replacement lemmas, in order to extend arguments involving mobile clusters.

Lastly, this work serves as a first effort in the direction of the larger question of classifying the types of diffusivities attainable through gradient dynamics. In a project currently in development, this is analysed in a constructive manner by identifying the collection of Bernstein models with a polynomial basis.

Finally, we note that an extension of the dynamics to long-range jumps cannot be performed with the same reasoning as for the PMM in \cite{LR:CDG2022}, since it breaks the gradient property. This is left as an open problem.

\section{Setting and main results}
\subsection{Context and main problem}

In the present section, we explain the context of this work and the questions addressed here at a more technical level. We start by fixing our setup. Let $N\in\mathbb{N}_+$ be fixed and consider the discrete lattice given by $\mathbb{T}_N:=\mathbb{T}/N\mathbb{Z}$, the one-dimensional discrete torus. The elements of $\mathbb{T}_N$ are referred to as "sites", and a pair of neighbouring sites, $\{x,y\}$ with $|x-y|=1$, as a "node". A configuration of particles is an element of the state space $\Omega_N=\{0,1\}^{\mathbb{T}_N}$ and will be recurrently denoted by $\eta$. We denote by $\eta(x)\in\{0,1\}$ the occupation value of $\eta\in\Omega_N$ at the site $x\in\mathbb{T}_N$. We now introduce the relevant operators associated with the process.  

\begin{Def}
In what follows, let $\eta\in\Omega_N$ , $f:\Omega_N\to\mathbb{R}$ and $x,y,z\in\mathbb{T}_N$ be arbitrary:
\begin{itemize}
        \item Let $\pi_x:\Omega_N\to\{0,1\}$ be the projection $\pi_x(\eta)=\eta(x)$;
        
        \item $\tau:\eta\mapsto\tau\eta$ is defined as the shift operator, $\pi_x(\tau\eta)=\pi_{x+1}(\eta)$, and short-write $\tau^i=\circ_{j=1}^i\tau$ for its $i$-th composition;

        \item Denote by $\theta_{x,y}$ the operator that exchanges the occupation-value of the sites $x,y$,
        \begin{align}
            (\theta_{x,y}\eta)(z)
            =\eta(z)1_{z\neq x,y}
            +\eta(y)1_{z=x}
            +\eta(x)1_{z=y};
        \end{align}

        \item For $\mathcal{O}$ any of the previously introduced operators, let $\mathcal{O} f(\eta):=f(\mathcal{O}\eta)$.
    \end{itemize}
    Moreover, introduce the linear operators $\nabla_{x,y},\nabla,\Delta$ through 
    \begin{align}
        \nabla_{x,y}[f]:=\theta_{x,y}f-f
        ,\quad
        \nabla[f]:=\tau f-f 
        \quad\text{and}\quad
        \Delta[f]=\nabla^2[\tau^{-1}f]
        .
    \end{align}
    Note that $\nabla$ corresponds to the forward difference operator and $\Delta$ to the discrete Laplacian operator.
        
    
\end{Def}
A space-homogenous exclusion process on $\mathbb{T}_N$ is characterized by its Markov generator, $\mathcal{L}$, defined through its action on $f:\Omega_N\to\mathbb{R}$ by
    \begin{align}
        \mathcal{L}_N[f](\eta)
        =
	\sum_{x\in\mathbb{T}_N}
        \mathbf{c}_{x,x+1}(\eta)
        \mathbf{e}_{x,x+1}(\eta)
        \nabla_{x,x+1}[f](\eta)
        +\sum_{x\in\mathbb{T}_N}
        \mathbf{c}_{x,x+1}(\eta)	
        \mathbf{e}_{x,x+1}(\eta)
        \nabla_{x,x+1}[f](\eta)
    \end{align}
where $\mathbf{c}_{x,x+1}\geq 0$ and $\mathbf{e}_{x,x+1}(\eta)=\eta(x)(1-\eta(x+1))$ are the \textit{kinetic} and \textit{exclusion} constraints for the hopping from the site $x$ to the site $x+1$;  and $\mathbf{c}_{x+1,x}$ and $\mathbf{e}_{x+1,x}$ the analogous constraints for a jump from $x+1$ to $x$. Both $\mathbf{c}_{x,x+1}\geq 0$ and $\mathbf{c}_{x+1,x}\geq 0$ depend on some range of interaction, possibly taking different values for different local configurations.

Our focus will be on symmetric dynamics in homogeneous media, that is, $\mathbf{c}_{x,x+1}=\mathbf{c}_{x+1,x}$ and $\mathbf{c}_{x,x+1}=\tau^x\mathbf{c}_{0,1}$. Conveniently, the generator is reversible with respect to the Bernoulli product measure $\nu_\rho^N$, given by the marginals $\nu_\rho^{N}(\eta\in\Omega_N:\;\eta(x)=1) = \rho$ for any $x\in\mathbb{T}_N$ and parametrised by $\rho\in[0,1]$.  

Expressing generically $\mathbf{c}_{0,1}=\tfrac12\mathbf{c}$ for an appropriate constraint $\mathbf{c}\geq 0$ allows us to write  
\begin{align}\label{generic:gen}  
    \mathcal{L}_N[f](\eta)  
    =  
    \sum_{x\in\mathbb{T}_N}  
    \mathbf{c}(\tau^x\eta)  
    \nabla_{x,x+1}[f](\eta),  
\end{align}  
associating univocally each symmetric process as just described with a constraint $\mathbf{c}$. Our goal is to introduce a specific constraint $\mathbf{c}$ with particular properties, filling a gap in the literature on the types of diffusion coefficients that can be derived from these discrete Markovian dynamics.

The key technical property of the forthcoming models is the so-called \textit{gradient property} \cite[II Subsection 2.4]{spohn:book}. In our context, this is reduced to the existence of some $\mathbf{H}:\Omega_N\to\mathbb{R}$ such that  
\begin{align}\label{generic:curr}  
    \mathbf{j}_{x,x+1}=-\nabla\tau^x\mathbf{H}  
    \quad\text{where}\quad  
    \mathbf{j}_{x,x+1}(\eta)=\mathbf{c}(\tau^x\eta)\mathbf{e}_{0,1}(\tau^x\eta)-\mathbf{c}(\tau^x\eta)\mathbf{e}_{1,0}(\tau^x\eta),  
\end{align}  
and $\mathbf{j}_{x,x+1}$ is the algebraic current associated with the node $\{x,x+1\}$. The current is, in fact, defined from the generator through $\mathcal{L}_N[\pi_x](\eta)=-\nabla\mathbf{j}_{x,x+1}(\tau^{-1}\eta)$, thus the gradient condition induces a "heat-like" equation $\mathcal{L}_N[\pi_x](\eta)=\Delta\mathbf{H}(\tau^x\eta)$. From Kolmogorov's equation, one then expects that the average path of a particle follows a discretised version of the evolution equation $\partial_t\rho=\partial_u^2\text{H}(\rho)$, where $\partial_t$ and $\partial_u$ correspond to the time and space derivatives, and $\text{H}$ coincides with the canonical ensemble, $\text{H}(\rho)=\int_{\Omega_N}\mathbf{H}\;\rmd\nu_\rho^N$. In particular, we also identify the diffusion coefficient $D(\rho):=\text{H}'(\rho)$ as $D(\rho)=\int_{\Omega_N}\mathbf{c}\;\rmd\nu_\rho^N$. This technical property is fundamental in different techniques for the study of the scaling limit of particle systems.




It is convenient to recall two well-known processes: the Symmetric Simple Exclusion Process (SSEP) \cite{SPITZ70,KL:book} and the (normalised) Porous Media Model (PMM) \cite{GLT}.  

\begin{Def}
    The SSEP is the process with generator given in \eqref{generic:gen} with $\mathbf{c}(\eta)=1$, for any $\eta\in\Omega_N$. We write its generator as $\mathcal{L}_N^{\text{SSEP}}$.

    For each fixed $n\in\mathbb{N}_+$, the PMM is the process given by the generator in \eqref{generic:gen} with $\mathbf{c}\equiv\mathbf{c}^n$, where
    \begin{align}\label{pmm-cons}
        \mathbf{c}^{n}(\eta):=\frac{1}{n+1}\sum_{i=0}^n
        \prod_{\substack{i=-n-1+j\\i\neq 0,1}}^j
        \eta(i)
        .
    \end{align}
In order to specify the parameter $n$, we shall write its generator as $\mathcal{L}_N^{n}$, and refer to the process by PMM($n$).
\end{Def}

The SSEP is a classical model for heat transfer, with a vast literature on it, and the PMM is a toy model for degenerate diffusion, enjoying relevant properties such as the existence of mobile clusters and blocked configurations \cite{GLT,BDGN}.  

From the previous observations, the PMM($n$) is associated with the diffusion coefficient $D_n(\rho)=\rho^n$. Moreover, whenever the constraint $\tilde{\mathbf{c}}:=a_0+a_1\mathbf{c}^1+a_2\mathbf{c}^2+\cdots$ is non-negative and conveniently bounded from above \cite{s2f,CGN24} (therefore under specific conditions on the coefficients $(a_i)_{i\geq0}$), a superposition principle holds, and $\tilde{\mathbf{c}}$ is associated with the diffusion coefficient $\tilde{D}(\rho)=a_0+a_1\rho+a_2\rho^2+\cdots$. This allows for a large variety of polynomials and functions that can be represented in series form to be associated with diffusion coefficients derived from gradient models. In particular, this provides a convenient microscopic context for the study of target properties of evolution equations.  

The technical gap that we identified in the classes of the aforementioned functions is that the PMMs do not allow the derivation of diffusion coefficients with two roots of multiplicity equal to or larger than two. Indeed, a polynomial in the standard basis $\{\rho^n\}_{n\geq 0}$ does not necessarily have non-negative coefficients, making it necessary to verify whether the underlying model is well-defined. Using the PMMs as a polynomial basis may lead to ill-defined models for surprising diffusivities, such as the diffusion coefficient $D_{n,L}(\rho)=\text{B}_{n,L}(\rho)$, indexed by $n,L\in\mathbb{N}$ with $n\leq L$ and given by  
\begin{align}  
    \text{B}_{n,L}(\rho)  
    =\binom{L}{n}\rho^n(1-\rho)^{L-n},  
\end{align}  
where $\{\text{B}_{n,L}(\rho)\}_{0\leq n\leq L}$ forms the Bernstein polynomial basis of degree $L$.  

The constraint $\tilde{\mathbf{c}}^{n-1,1}:=\mathbf{c}^{n-1} - \mathbf{c}^{n}$ is well-defined and associated with $\rho^{n-1}(1-\rho)$, allowing for a simple root and another of arbitrary multiplicity. This is a consequence of the sequence $(\mathbf{c}^{n})_{n\geq 0}$ being non-increasing \cite[Proposition 2.18]{s2f}. However, the sequence $(\tilde{\mathbf{c}}^{n,1})_{n\geq 0}$ is not monotone, and in general, for any $k\geq 2$, there exist configurations such that  
\begin{align}  
    \tilde{\mathbf{c}}^{n,k}(\eta):=\sum_{j=0}^k\binom{k}{j}(-1)^n\mathbf{c}^{n+j}(\eta)<0.  
\end{align}  
An example is $\eta$ where every site is occupied except for the site $-n-k$, corresponding to $\tilde{\mathbf{c}}^{n,k}(\eta)=-(n+k+1)^{-1}$.  

Concluding the present discussion, we address in this manuscript two questions: Is there a gradient dynamics associated with the diffusion coefficient $D_{n,L}(\rho)$? And if so, what is the dynamics associated with the underlying polynomial basis $\{\rho^n\}_{n\geq 0}$?

\subsection{Main Results and Outline of the paper}

    
In order to present the models we need to introduce some relevant definitions.  
    
\begin{Def}
    In what follows, fix $L\in\mathbb{N}_+$, let $0\leq j\leq L$ and $\eta\in\Omega_N$ be arbitrary:
    \begin{itemize}
        \item Introduce the set $W_j^L:=\llbracket -j,-j+L+1\rrbracket\backslash\{0,1\}$; 
        \item Let $m_j^L(\eta)$ be the number of particles in $W_j^L$;

        \item Let $\inner{\eta}_{W_j^L}=\tfrac{1}{L}m_j^L(\eta)$ be the density of particles in the box $W_j^L$.
    \end{itemize}
\end{Def}

We are now ready to introduce the main models.

\begin{Def}[Bernstein model]\label{def:B}
    For each $ n,L\in\mathbb{N}$ fixed with $n\leq L$, define the constraint $\mathbf{b}_{n,L}:\Omega_N\to[0,1]$ as 
    \begin{align}\label{rate:b}
		\mathbf{b}_{n,L}
		:=
		\frac{1}{L+1}\sum_{j=0}^{L}
		\mathbf{b}_{n,L}^j
        \quad\text{with}\quad
            \mathbf{b}_{n,L}^j(\eta)
            :=\mathbf{1}\big\{\inner{\eta}_{W_j^L}=\tfrac{n}{L}\big\}.
    \end{align}
    We refer to the process induced by setting $\mathbf{c}=\mathbf{b}_{n,L}$ in \eqref{generic:gen} as B($n$,$L$), and its corresponding Markov generator as $\mathcal{L}_N^{n,L}$.

\end{Def}

In Figure \ref{fig:PMMnk} we present an example of the transition rates. In general, assuming that $\eta(0)+\eta(1)=1$, the rate $\mathbf{b}_{n,L}(\eta)(\mathbf{e}_{0,1}(\eta)+\mathbf{e}_{1,0}(\eta)$ coincides with the proportion of boxes of length $L+2$ containing the node $\{0,1\}$ with exactly $n+1$ particles.

We now introduce the reduced PMM.
\begin{Def}[Reduced Porous Media Model]\label{def:P}
    For each fixed natural numbers $\ell \leq L<N$, each $0\leq j\leq L$ and each $\eta\in\Omega_N$, let 
\begin{align}\label{new-base}
    \mathbf{p}_{\ell;L}
    &:=
    \frac{1}{L+1}
    \sum_{j=0}^{L}
    \mathbf{p}_{\ell;L}^j
    \quad\text{with}\quad
    \mathbf{p}_{\ell;L}^j(\eta)
    :=
    \frac{\binom{m_j^L(\eta)}{\ell}}{\binom{L}{\ell}}
    \mathbf{1}
    \{
    m_j^L(\eta)\geq \ell
    \}
    .    
\end{align}
    We refer to the process induced by setting $\mathbf{c}=\mathbf{p}_{\ell;L}$ in \eqref{generic:gen} as $\text{PMM}_{L}(\ell)$, and its corresponding Markov generator as $\mathcal{L}_N^{\ell;L}$.
\end{Def}
Note that for each $0\leq j\leq L$ the constraint $\mathbf{p}_{\ell;L}^j$ corresponds to the probability of choosing $\ell$ out of the $m_j^L$ particles in the box $W_j^L$. We allude the reader to Figure \ref{fig:redPMM}.

The main processes are related to the PMM through B($n$,$n$) and $\text{PMM}_n(n)$ being identified with the PMM($n$). Moreover, it is simple to identify the diffusion coefficients $\int_{\Omega_N}\mathbf{B}_{n,L}\rmd\nu_\rho^N=\text{B}_{n,L}(\rho)$ and $\int_{\Omega_N}\mathbf{p}_{\ell;L}\rmd\nu_\rho^N=\rho^\ell$ (see Proposition \ref{prop:P-prop} for the latter). In particular, fixed $L$, the $\text{PMM}_L(\ell)$ corresponds to a reduction of the diffusivity of the PMM($L$),  $\alpha^L\to\alpha^\ell$, while maintaining the interaction range $2L+2$.

Section \ref{sec:char} is devoted to relating the two previous models, showing the gradient property and other elementary properties. Specifically, in subsection \ref{subsec:bin} we show that the models are related through a monomial decomposition and inversion formula. The next Proposition summarizes their connections in this work.
\begin{Prop}\label{prop:identity}
Fixed $L\in\mathbb{N}$, it holds that 
    \begin{align}
    \mathbf{b}_{n,L}
    &=
    \sum_{\ell=n}^{L}
    (-1)^{\ell-n}
    \binom{L}{\ell}\binom{\ell}{n}
    \mathbf{p}_{\ell;L}
    &&\text{for each } 0\leq n\leq L
    ,\label{b-decomp}
    \\
    \mathbf{p}_{\ell;L}
    &=\sum_{n=\ell}^L
    \frac{\binom{n}{\ell}}{\binom{L}{\ell}}
    \mathbf{b}_{n,L},
    &&\text{for each } 0\leq \ell\leq L
    \label{p-inv}
    \\
    \mathbf{b}_{n,L}(\eta)
    &\leq
    \binom{L}{n}
    \mathbf{p}_{n;L}(\eta),
    &&\text{for each } \eta\in\Omega_N
    .
    \label{ineq}
    \end{align}
\end{Prop}
The identification \eqref{p-inv} is very simple and explained in Lemma \ref{lem:inverse}. The decomposition \eqref{b-decomp} is not obvious by an application of the inverse binomial transformation because each term of $(\mathbf{p}_{\ell,L})_{\leq L}$ depends on the choice of $L$ a priori. This is shown in Lemma \ref{lem:decomp}, but it can also be derived straightforwardly from the principle of inclusion-exclusion as presented in \cite[Equation (4.2.3)]{book:genfun}. The inequality \eqref{ineq} is a simple observation directly from the definition of the constraints.

A key result is the verification of the gradient property of the Bernstein model, that we state below and prove in subsection \ref{subsec:grad}
\begin{Prop}\label{prop:grad-B}
Fixed $L\in\mathbb{N}$, shorten $\inner{\eta}_L=\tfrac{1}{L+1}\sum_{z\in \llbracket 0,L\rrbracket}\eta(z)$ for the density in the box $\llbracket 0,L\rrbracket$. For any $n\leq L$ it holds that
\begin{align}
    \mathbf{b}_{n,L}(\eta)
    (\mathbf{e}_{0,1}(\eta)-\mathbf{e}_{1,0}(\eta))
    =-\nabla\mathbf{H}_{n,L}(\eta)
\end{align}
with $\mathbf{H}_{n,L}=\mathbf{h}_{n,L}+\mathbf{g}_{n,L}$ and
\begin{align}\label{grad:c}
    \mathbf{h}_{n,L}(\eta)
    &=
    \frac{1}{L+1}\mathbf{1}{
    \big\{
    \inner{\eta}_{L}\geq\tfrac{n+1}{L+1}
    \big\}
    }
    %
    ,
    \\
    \mathbf{g}_{n,L}(\eta)
    &=\frac{1}{L+1}\sum_{j=1}^{L}
    \sum_{i=1}^{j}
    \mathbf{b}_{n,L}^j(\tau^{-i}\eta)
    (\mathbf{e}_{0,1}(\tau^{-i}\eta)
    -\mathbf{e}_{1,0}(\tau^{-i}\eta)
    )
    .
\end{align}
\end{Prop}

Provided the proposition above and the relation of $\mathbf{p}_{\ell;L}$ with the inverse binomial transformation of the sequence $(\mathbf{b}_{n,L})_{n\leq L}$ (see \eqref{p-inv}), it is clear that the $\text{PMM}_L(\ell)$ is a gradient model as-well. The precise expressions for the quantities involved in the gradient property, for the reduced PMM, can be identified from a decomposition of the ones associated with the Bernstein model. This is performed in Lemma \ref{lem:p-grad}. 

In subsection \ref{subsec:further} we prove other properties of the main models, both dynamic and algebraic. Concretely, we prove Propositions \ref{prop:B-prop} and \ref{prop:P-prop} below, summarizing also all the other more elementary properties of these models.

\begin{Prop}\label{prop:B-prop}
For each $n,L\in\mathbb{N}$ with $L\geq n$, the B($n$,$L$) enjoys the following properties:
    \begin{enumerate}[label=(\roman*)]

    \item \textit{Partition of the unity}: for any integer $ 0<L<N/2 $ it holds that
    \begin{align}
    \mathcal{L}_N^{\text{SSEP}}
    =\sum_{n=0}^L\mathcal{L}_N^{n,L};
    \end{align}

    \item \textit{Symmetry}: for any $\eta\in \Omega_N$ it holds that $\mathbf{b}_{n,L}(\eta)=\mathbf{b}_{L-n,L}(\overline{\eta})$, with $\eta\mapsto\overline{\eta}$ the map that exchanges the particles by vacant sites and vice-versa;

    \item \textit{Interpolation}: for each $n\in\mathbb{N}_+$ it holds that \textnormal{B($n$,$n$)$=$PMM($n$)};

    \item \textit{Range}: For each $\eta\in\Omega_N$ it holds that $\mathbf{b}_{n,L}(\eta)\in [0,1]$. 

    \item \textit{Blocked configurations}: there is $\eta\in\Omega_N$ such that $\mathbf{b}_{n,\ell}(\tau^x\eta)=0$ for any $x\in\mathbb{T}$;
    
    \item \textit{Mobile Clusters}: any box of length $L+2$ with exactly $n+1$ particles constitutes a mobile cluster. 
\end{enumerate}
\end{Prop}

\begin{Prop}\label{prop:P-prop}
    For each natural number $L<N$ fixed and any $0\leq \ell\leq L$, it holds that 
\begin{enumerate}[label=(\roman*)]
        \item \textit{Mobile Clusters}: any box of length $L+2$ with \textit{at least} $\ell+1$ particles and at least one vacant site is a mobile cluster;
        
        \item There are blocked configurations; 

        \item Interpolation: $\text{PMM}_L(L)=\text{PMM}(L)$;

        \item \textit{Range}: For each $\eta\in\Omega_N$ it holds that $\mathbf{p}_{\ell;L}(\eta)\in [0,1]$. 
        
        \item Monotony: the sequence $(\mathbf{p}_{\ell;L})_{0\leq \ell\leq L}$ is non-increasing;
        
        \item Equilibrium expected value: for any $\rho\in[0,1]$, it holds that  $\int_{\Omega_N}
            \mathbf{p}_{\ell;L}\rmd\nu_\rho^N=\alpha^\ell$.
    \end{enumerate}
\end{Prop}

\section{Characterization of the collection of models}\label{sec:char}

\subsection{Binomial Transformation}\label{subsec:bin}

We start by showing Proposition \ref{prop:identity}, relating the Bernstein and reduced PMM via the binomial transformation.

\begin{Lemma}[Inversion formula]\label{lem:inverse}
For each $n,L\in\mathbb{N}$ with $L>n$, it holds that
    \begin{align}\label{p-inverse}
    \mathbf{p}_{\ell;L}(\eta)
    &=\sum_{n=\ell}^L
    \frac{\binom{n}{\ell}}{\binom{L}{\ell}}
    \mathbf{b}_{n,L}.
\end{align}
\end{Lemma}
\begin{proof}
    It is enough to express, for each $0\leq j\leq L$,
\begin{align}
    \mathbf{p}_{\ell;L}^j(\eta)
    =    
    \frac{\binom{m_j^L(\eta)}{\ell}}{\binom{L}{\ell}}
    \mathbf{1}_{
    \{
    m_j^L(\eta)\geq \ell
    \}
    }
    =\sum_{n=\ell}^L
    \frac{\binom{n}{\ell}}{\binom{L}{\ell}}
    \mathbf{1}_{
    \{
    m_j^L(\eta)=n
    \}
    }
\end{align}
and identify $\mathbf{1}_{
    \{
    m_j^L(\eta)=n
    \}
    }
    =\mathbf{b}_{n,L}^j(\eta)$. 
\end{proof}

In order to prove the next lemma it will be convenient to introduce some notation.
\begin{Def}\label{def:sets}
    Let $A,P\subseteq\mathbb{T}_N$ , $a\in\mathbb{T}_N$ and $L\in\mathbb{N}$ be arbitrarily fixed. We write ${\eta^P(a)}:=\mathbf{1}_{a\notin P}\eta(a)+\mathbf{1}_{a\in P}(1-\eta(a))$ and ${\eta(A)}:=\prod_{a\in A}\eta(a)$. Consequently, $\eta^P(A)=\prod_{a\in A}\eta^P(a)$. Moreover, let {$\mathcal{P}_L(A)$} denote all the subsets of $A$ with exactly $L$ elements. For $L>|A|$ and $L=0$ by convention $\mathcal{P}_L(A)=\emptyset$;


\end{Def}

\begin{Lemma}[Transformation to monomials]\label{lem:decomp}
For each $n,L\in\mathbb{N}$ with $L>n$, it holds that
    \begin{align}\label{b-p_base}
\mathbf{b}_{n,L}=
\sum_{\ell=n}^{L}
(-1)^{\ell-n}
\binom{L}{\ell}\binom{\ell}{n}
\mathbf{p}_{\ell;L}
.
\end{align}
\end{Lemma}

\begin{proof}
Recalling \eqref{rate:b}, for each $0\leq j\leq L$, from the distributive rule one can express 
\begin{align}\label{cons:monomial}	
    \mathbf{b}_{n,L}^j
    &=
    \sum_{P\in \mathcal{P}_{L-n}(W_j^L)}
    \eta^{P}(W_j^L)
    \\
    &=
    \sum_{P\in \mathcal{P}_{L-n}(W_j^L)}
    \sum_{\ell=0}^{L-n}(-1)^\ell
    \sum_{Q\in\mathcal{P}_{\ell}(P)}
    \eta([W_j^L\backslash P]\cup Q )
    \\
    &=
    \sum_{\ell=0}^{L-n}
    (-1)^\ell
    \binom{n+\ell}{\ell}
    \sum_{P\in\mathcal{P}_{n+\ell}(W_j)}\eta(P).
\end{align} 
The third equality is consequence of, fixed $\ell\leq L, P\in\mathcal{P}_{L-n}(W_j^L)$ and $Q\in\mathcal{P}_{\ell}(P)$, clearly $[W_j^L\backslash P]\cup Q\in \mathcal{P}_{n+\ell}(W_j^L)$; and for each $P_\ell\in\mathcal{P}_{n+\ell}(W_j^L)$ there are $\binom{n+\ell}{\ell}$ pairs $(P,Q_P)$ such that $P\in P_{L-n}(W_j^L)$ and $Q_P\in\mathcal{P}_{\ell}(P)$, where $P_\ell=[W_j^L\backslash P]\cup Q_P$.

One sees that for each $0\leq\ell\leq L-n$,
    \begin{align}
        \sum_{P\in\mathcal{P}_{n+\ell}(W_j^L)}
	\eta(P)
        =
        \binom{m_j^L(\eta)}{n+\ell}
        \mathbf{1}_{
        \{
        m_j^L(\eta)
        \geq
        n+\ell 
        \}
        },    
    \end{align}
leading to
\begin{align}\label{bj-mono}
    \mathbf{b}_{n,L}^j(\eta)
=    \sum_{\ell=0}^{L-n}(-1)^\ell
\binom{n+\ell}{\ell}\binom{L}{n+\ell}
    \frac{\binom{m_j^L(\eta)}{n+\ell}}{\binom{L}{n+\ell}}
        \mathbf{1}_{
        \{
        m_j^L(\eta)
        \geq
        n+\ell 
        \}
        }
        .
\end{align}
Summing the above over $0\leq j\leq L$ and performing a change of variables we identify the collection of maps $(\mathbf{p}_{\ell;L})_{0\leq \ell\leq L-n}$ as in \eqref{new-base}
\end{proof}

\subsection{Gradient condition}\label{subsec:grad}
We now verify the gradient property for the B($n$,$L$).
\begin{proof}[Proof of Proposition \ref{prop:grad-B} ]

We recall Definition \ref{def:sets} and the constraint as in \eqref{rate:b}. Let us introduce the set $M_L=\llbracket0,L+1\rrbracket$, and note that  
    \begin{align}
        \mathbf{b}_{n,L}(\eta)
        =\frac{1}{L+1}
        \sum_{j=0}^{L}
        \mathbf{s}_{n,L}^j(\tau^{-j}\eta)
        \quad\text{where}\quad
        \mathbf{s}_{n,L}^j(\eta)
        =\sum_{\substack{P\in \mathcal{P}_{L-n}(M_L)\\j,j+1\notin P}}
        \eta^{P}(M_L\backslash\{j,j+1\}).
    \end{align}
    Expressing $\mathbf{1}-\tau^{-j}=\sum_{i=1}^{j}\nabla\tau^{-i}$, we can rewrite
    \begin{align}\label{grad:step1}
        \mathbf{b}_{n,L}(\eta)(\mathbf{e}_{1,0}(\eta)-\mathbf{e}_{0,1}(\eta))
        &=
        \frac{1}{L+1}\sum_{j=0}^{L}
    \sum_{\substack{P\in \mathcal{P}_{L-n}(M_L)\\j,j+1\notin P}}
\big\{
    \eta^{P\cup\{j+1\}}(
        M_L
        )
        -
        \eta^{P\cup\{j\}}(
        M_L
        )
\big\}
    \\&+\frac{1}{L+1}\nabla\mathbf{g}^{n,L}(\eta),
    \end{align}
with $\mathbf{g}^{n,L}$ as in \eqref{grad:c}. Let us introduce the auxiliary sets
\begin{align}
    \mathcal{W}_P
    &=\{x:x\notin P, x+1\in P,\;x,x+1\in M_L\}
    ,\\
    \Tilde{\mathcal{W}}_P
    &=\{x:\;x\in P,x+1\notin P,\;x,x+1\in M_L\}.
\end{align}
We claim that the first term in the right-hand side of \eqref{grad:step1} can be further expressed as 
\begin{multline}\label{grad:step2}
    \frac{1}{L+1}
    \sum_{\substack{P\in \mathcal{P}_{L-n+1}(M_L)\\P\neq \{0,\dots,L-n\}}}
    |\mathcal{W}_P|\eta^{P}(M_L)
    -\sum_{\substack{P\in \mathcal{P}_{L-n+1}(M_L)\\P\neq \{L+1-(L-n),\dots,L+1\}}}
    |\Tilde{\mathcal{W}}_P|\eta^{P}(M_L)
\\=    
    \frac{1}{L+1}
    \sum_{P\in \mathcal{P}_{L-n+1}(M_L)}
    |\mathcal{W}_P|\eta^{P}(M_L)
    -\sum_{P\in \mathcal{P}_{L-n+1}(M_L)}
    |\Tilde{\mathcal{W}}_P|\eta^{P}(M_L).
\end{multline}
In order to show this, first note that 
\begin{align}
    \bigcup_{j=0}^L 
    \big\{
    P\cup \{j\}
    \big\}_{P\in \mathcal{P}_{L-n}(M_L\backslash\{j,j+1\})}
    &=\mathcal{P}_{L-n+1}(M_L)\backslash\big\{0,\dots,L-n \big\},
    \\
    \bigcup_{j=0}^L 
    \big\{
    P\cup \{j+1\}
    \big\}_{P\in \mathcal{P}_{L-n}(M_L\backslash\{j,j+1\})}
    &=\mathcal{P}_{L-n+1}(M_L)\backslash\big\{L+1-(L-n),\dots,L+1 \big\}
    .
\end{align}
Next, note that fixed $Q\in \mathcal{P}_{L-n+1}(M_L)\backslash\{0,\dots,n\}$, for each $j'\in W_Q$ there is a unique $P_{j'}\in \mathcal{P}_{L-n}(M_L\backslash\{j',j'+1\})$ such that we can express 
\begin{align}
    Q=P_{j'}\backslash\{j',j'+1\}\cup\{j'+1\}.
\end{align}
In particular, there is a total of $|\mathcal{W}_Q|$ pairs $(j',P_{j'})$ as above. This proves the first summation in \eqref{grad:step2}. For the second summation, the reasoning is analogous and one should consider the set $\Tilde{W}_Q$ instead. This proves the claim.

Established \eqref{grad:step2}, it is straightforward to see that for any $P\in \mathcal{P}_{L-n+1}(M_L)$, 
\begin{align}\label{grad:step3}
    |\mathcal{W}_P|-|\Tilde{\mathcal{W}}_P|
    =
    \begin{cases}
    1, &0\notin P,\;L+1\in P,\\
    -1, &0\in P,\;L+1\notin P,\\
    0, &\text{otherwise}.
    \end{cases}
\end{align}
Indeed, let us map $P\in \mathcal{P}_{L-n+1}(M_L)\mapsto \Omega_L=\{0,1\}^{M_L}$ with the identification $\xi(p)=1\Leftrightarrow p\in P$ and $\xi(p)=0\Leftrightarrow p\notin P$. Then $|\mathcal{W}_P|=\{\{x,x+1\}\subset M_L:\; \xi(x)=1,\;\xi(x+1)=0\}$ and, analogously, $|\Tilde{M}_p|=\{\{x,x+1\}\subset M_L:\; \xi(x)=0,\;\xi(x+1)=1\}$. Each local configuration of the form $(1,0,\dots,0,1)$ within $\xi$, with (say) $h>0$ successive vacant sites contributes to $+1$ for both $|\mathcal{W}_P|$ and $|\Tilde{\mathcal{W}}_P|$. If $h=0$, the local configuration does not contribute. The same observation is valid for $(0,1,\dots,1,0)$ with $q>0$ or $q=0$ consecutive particles. In this manner, the values of $|\mathcal{W}_P|$ and $|\Tilde{\mathcal{W}}_P|$ may differ depending on the occupations at the boundary sites $0$ and $L+1$. From here one notes that if $\xi(0)=1$ then $|\mathcal{W}_P|=|\Tilde{\mathcal{W}}_P|$ if $\xi(L+1)=1$; or $|\mathcal{W}_P|=|\Tilde{\mathcal{W}}_P|+1$ if $\xi(L+1)=0$. Likewise, if $\xi(0)=0$ then $|\Tilde{\mathcal{W}}_P|=|\mathcal{W}_P|$ if $\xi(L+1)=0$; or $|\Tilde{\mathcal{W}}_P|=|\mathcal{W}_P|+1$ if $\xi(L+1)=1$.    

In this way, from \eqref{grad:step3} we can express further \eqref{grad:step2} as 
\begin{align}
    \sum_{\substack{
    P\in \mathcal{P}_{L-n+1}(M_L)
    }}
    (
    \mathcal{W}_P-\Tilde{\mathcal{W}}_P
    )
    \eta^{P}(M_L)
    &=
    \sum_{\substack{
    P\in \mathcal{P}_{L-n+1}(M_L)
    \\0\notin P, L+1\in P 
    }}
    \eta^{P}(M_L)
    \\&-
    \sum_{\substack{
    P\in \mathcal{P}_{L-n+1}(M_L)
    \\0\in P, L+1\notin P 
    }}
    \eta^{P}(M_L)
\end{align}

Let us write the right-hand side above as $\alpha_{L-n+1}(\eta)$. The next step is to  identify the relation
\begin{align}\label{rec:rel}
    \begin{cases}
        \alpha_{i+1}=\alpha_{i}-\nabla\beta_{i}
        ,& 1\leq i \leq L-n\\
    \alpha_{1}(\eta)=-\nabla\eta(M_{L-n}),&
    \end{cases}
\quad\text{with}\quad
    \beta_{i}
    =\sum_{
    P\in\mathcal{P}_{i}(M_{n+i-2})
    }
    \eta^{P}(M_{n+i-1})
    .
\end{align}
Indeed, shortening $k\equiv L-n$,
\begin{align}
    \alpha_{k+1}(\eta)
    &=\sum_{\substack{
    P\in\mathcal{P}_{k}(M_{L-1})\\0\notin P
    }}
    \eta^P(M_{L-1})(1-\eta(L+1))
    -
    \sum_{\substack{
    P\in\mathcal{P}_{k}(M_{L-1})\\0\notin P
    }}
    (1-\eta(0))\eta^P(M_{L}\backslash\{0\})
    \\&
    =\sum_{\substack{
    P\in\mathcal{P}_{k}(M_{L-1})\\0\notin P
    }}
    \left\{
    \eta^P(M_{L-1})-\eta^P(M_L\backslash\{0\})
    \right\}
    ,
\end{align}
and because
\begin{align}
    \sum_{\substack{
    P\in\mathcal{P}_{k}(M_{L-1})
    \\0\notin P
    }}
    \eta^P(M_L\backslash\{0\})
    =
    \sum_{
    P\in\mathcal{P}_{k}(M_{L-2})
    }
    \eta^{1+P}(1+M_{L-1})
    =\beta_{k}(\tau^1\eta)
    ,
\end{align}
summing and subtracting $\beta_{k}(\eta)$, one obtains
\begin{align}
    \alpha_{k+1}(\eta)
    &=
    \sum_{\substack{
    P\in\mathcal{P}_{k}(M_{L-1})
    \\0\notin P
    }}
    \eta^P(M_{L-1})
    -
    \sum_{
    \substack{
    P\in\mathcal{P}_{k}(M_{L-2})\\L\notin P
    }
    }
    \eta^{P}(M_{L-1})
    -\nabla\beta_{k}(\eta)
    \\&
    =\alpha_{k}(\eta)-\nabla\beta_{k}(\eta)
    ,
\end{align}
with the last equality justified from the terms in the first summation where $L\in P$ cancelling with the terms in the second summation where $0\in P$.

Up to this point, from \eqref{rec:rel} we have that $\mathbf{b}_{n,L}(\eta)(\mathbf{e}_{1,0}(\eta)-\mathbf{e}_{0,1}(\eta))=\nabla\mathbf{H}_{n,L}(\eta)$, with $\mathbf{H}_{n,L}=\mathbf{h}_{n,L}+\mathbf{g}_{n,L}$ and  
    \begin{align}\label{h:decomp}
    \mathbf{h}_{n,L}(\eta)
    =
    \frac{1}{L+1}\sum_{i=0}^{L-n}    
    \sum_{
    P\in\mathcal{P}_{i}(\llbracket0,n+i-1\rrbracket)
    }
    \eta^{P}(\llbracket0,n+i\rrbracket)
    .
    \end{align}
We then see that
    \begin{align}
    \sum_{i=0}^{L-n}    
    \sum_{
    P\in\mathcal{P}_{i}(\llbracket0,n+i-1\rrbracket)
    }
    \eta^{P}(\llbracket0,n+i\rrbracket)
    &=
    \sum_{i=0}^{L-n}
    \eta(n+i)
    \sum_{
    P\in\mathcal{P}_{i}(\llbracket0,n+i-1\rrbracket)
    }
    \eta^{P}(\llbracket0,n+i-1\rrbracket)
    \\
    &=
    \sum_{i=0}^{L-n}
    \mathbf{1}{
    \big\{
    \inner{\eta}_{n+i-1}=\tfrac{n}{n+i}
    \big\}
    }
    \eta(n+i)
    \\
    &=
    \mathbf{1}{
    \big\{
    \inner{\eta}_{L}\geq \tfrac{n+1}{L+1}
    \big\}
    }
    ,
\end{align}
concluding the proof. 
\end{proof}

The gradient property for the $\text{PMM}_L$ is consequence of the decomposition \eqref{b-p_base} and the gradient property of B($n$,$L$). However, these relations do not provide a simple explicit expression for the related quantities. This is the content of the next lemma. 
\begin{Lemma}\label{lem:p-grad}
For each $n,L\in\mathbb{N}$ with $L\geq n$ it holds that
\begin{align}
    \mathbf{p}_{\ell;L}(\eta)
    (\mathbf{e}_{0,1}(\eta)-\mathbf{e}_{1,0}(\eta))
    =-\nabla\mathbf{H}_{\ell;L}(\eta)
\end{align}
with $\mathbf{H}_{\ell;L}=\mathbf{h}_{\ell;L}+\mathbf{g}_{\ell;L}$ and
\begin{align}
    \mathbf{h}_{\ell;L}(\eta)
    &
    =
    \frac{1}{L+1}\frac{1}{\binom{L}{\ell}}
    \mathbf{1}_{
    \{
    \inner{\eta}_{L}\geq \tfrac{\ell+1}{L+1}
    \}
    }    ,
    \\
    \mathbf{g}_{\ell;L}(\eta)
    &=\frac{1}{L+1}\sum_{j=1}^{L}
    \sum_{i=1}^{j}
    \mathbf{p}_{\ell;L}^j(\tau^{-i}\eta)
    (\mathbf{e}_{0,1}(\tau^{-i}\eta)
    -\mathbf{e}_{1,0}(\tau^{-i}\eta)
    )
    .
\end{align}
\end{Lemma}
\begin{proof}
Recalling Proposition \ref{prop:grad-B} and the decomposition \eqref{bj-mono}, we have that 
\begin{align}
    \mathbf{g}_{n,L}(\eta)
    &=\sum_{\ell=n}^L(-1)^{\ell-n}\binom{L}{n}\binom{L-n}{\ell-n}
    \frac{1}{L+1}\sum_{j=1}^{L}
    \sum_{i=1}^{j}
    \mathbf{p}_{\ell;L}^j(\tau^{-i}\eta)
    (\mathbf{e}_{0,1}(\tau^{-i}\eta)
    -\mathbf{e}_{1,0}(\tau^{-i}\eta)
    )
    ,
\end{align}
being then enough to decompose $\mathbf{h}_{n,L}$. Precisely, because $\binom{L}{n}\binom{L-n}{\ell-n}=\binom{L}{\ell}\binom{\ell}{n}$, this corresponds to showing the identity 
\begin{align}\label{connect}
    \mathbf{1}_{
    \big\{
    \inner{\eta}_{L}\geq\tfrac{n+1}{L+1}
    \big\}
    }
    =
    \sum_{\ell=n}^L
    (-1)^{\ell-n}
    \binom{\ell}{n}
    \mathbf{1}_{
    \{
    \inner{\eta}_{L}\geq \tfrac{\ell+1}{L+1}
    \}
    },
\end{align}
where we note that the left-hand side above equals $(L+1)\mathbf{h}_{n,L}(\eta)$. This can be seen as a direct application of the principle of inclusion-exclusion, or equivalently, the quantity in the left-hand side above equals
    \begin{align}\label{decomp-h}
\sum_{i=n}^{L}    
    \sum_{
    P\in\mathcal{P}_{i-n}(\llbracket0,i-1\rrbracket)
    }
    \eta^{P}(\llbracket0,i\rrbracket)
    &=
    \sum_{i=n}^{L}    
    \eta(i)
    \sum_{\ell=0}^{i-n}
    (-1)^{\ell}
    \sum_{
    \substack{
    P\in\mathcal{P}_{i-n}(\llbracket0,i-1\rrbracket)
    \\
    Q\in\mathcal{P}_{\ell}(P)
    }}
    \eta([\llbracket0,i\rrbracket\backslash P]\cup Q)
    \\
    &=
    \sum_{i=n}^{L}    
    \sum_{\ell=n}^{i}
    (-1)^{\ell-n}
    \binom{\ell}{n}
    \sum_{
    P\in\mathcal{P}_{\ell}(\llbracket0,i-1\rrbracket)
    }
    \eta(P)\eta(i)
    .
    \end{align}
Inverting the summations over $i$ and $\ell$ above and identifying 
\begin{align}
    \sum_{i=\ell}^L
    \sum_{
    P\in\mathcal{P}_{\ell}(\llbracket0,i-1\rrbracket)
    }
    \eta(P)\eta(i)
    =
    \mathbf{1}_{
    \{
    \inner{\eta}_{L}\geq \tfrac{\ell+1}{L+1}
    \}
    }
\end{align}
concludes the proof.
\end{proof}

\subsection{Further properties}\label{subsec:further}

In this subsection we prove Propositions \ref{prop:B-prop} and \ref{prop:P-prop}. The former is related to the Bernstein model, while the latter with the reduced PMM. In order to do so, let us recall the notion of mobile cluster \cite{GLT} used here. 
\begin{Def}\label{def:mc}
Let $\blacksquare$ represent a fixed finite box, composed by particles and vacant sites. Representing a particle by $\bullet$ and a hole by $\circ$, we say that $\blacksquare$ constitutes a mobile cluster if it enjoys the following.
\begin{itemize}[label=$\triangleright$]
\item \textit{Mobility}: the transitions $\blacksquare\bullet \leftrightarrow \bullet\blacksquare$ and $\blacksquare\circ\leftrightarrow\circ\blacksquare$ are possible with a finite number of jumps, and independent of the rest of the configuration; 
    \item \textit{Mass transport}: it is always possible for a jump to occur in a node, if there exists a cluster in the vicinity of the respective node, that is, $\blacksquare\circ\bullet\leftrightarrow\blacksquare\bullet\circ$ and $\circ\bullet\blacksquare\leftrightarrow\bullet\circ\blacksquare$.
\end{itemize}	    
\end{Def}







    

\begin{proof}[Proof of Proposition \ref{prop:B-prop}]
    For \textit{(i)}, simply note that for any $ \eta\in\Omega_N $ fixed and any $0\leq j\leq L$,
    \begin{align}
    \sum_{n=0}^L\mathbf{b}_{n,L}^j(\eta)
    =
    \sum_{n=0}^L\mathbf{1}\big\{\inner{\eta}_{W_j}=\tfrac{n}{L}\big\}
    =1
    .
    \end{align}
Property \textit{(ii)} is a simple observation that $n$ particles corresponds to $L-n$ vacant sites. The interpolation, \textit{(iii)} follows from $\mathbf{b}_{n,n}=\mathbf{c}^{n}$, with the latter as in \eqref{pmm-cons}, by direct substitution of $L$ by $n$ in \eqref{rate:b} and identifying $\mathbf{1}\{\inner{\eta}_{W_j}=1\}=\eta(W_j)$.

On \textit{(iv)}, clearly $\mathbf{b}_{n,L}(\eta)\in[0,1)$ for any $\eta\in\Omega_N$, and we just need to observe that the constraint man attain the maximum value. Provided any box of length $L$ with exactly $n$ particles, one can construct $\xi$ such that $\mathbf{b}_{n,L}(\xi)=1$ by defining $\xi$ restricted to $W_0^L$ equal to the configuration in the aforementioned box; then defining recursively the configuration in $W_1^L,W_2^L,\cdots$ so that each of these windows contain exactly $n$ particles.

Regarding \textit{(iv)}, an example is the configuration where each particle is at a distance larger than $ L $ from each other; or each vacant site is at a distance larger than $ L $ from each other.

On $(vi)$, recall Definition \ref{def:mc} and let $\Box$ represent a cluster of particles composed by a box of length $L$ with exactly $n$ particles, let $\bullet$ represent a particle and $\circ$ a hole, and let $\blacksquare$ represent a box of length $L+2$ with exactly $n+1$ particles. We note that the particles in $\blacksquare$ can be reorganized arbitrarily in a finite number of steps: any pair $\circ\bullet$ (resp. $\bullet\circ$) inside $\blacksquare$ can transition to $\bullet\circ$ (resp. $\circ\bullet$) because it is contained inside a box with $n+1$ particles. As a consequence of this observation, both the mobility and mass transport properties directly follow by convenient reorganizations of the cluster.
\end{proof}

We conclude this work by showing Proposition \ref{prop:P-prop}.



        
        

\begin{proof}[Proof of Proposition \ref{prop:P-prop}]
On (i), from the formula \eqref{p-inverse}, the process corresponds to a superposition of weighted Bernstein models, therefore any mobile cluster of any of the processes $\{$B($n$,$L$)$\}_{\ell\leq n\leq L}$ constitutes also a mobile cluster for the $\text{PMM}_L(\ell)$. In particular, by union, mobile clusters are composed by boxes of length $L+2$ with \textit{at least} $\ell+1$ particles and at least one vacant site. 

The existence of blocked configurations, (ii), is also clear and any configuration composed by boxes of length $L$, with less than $\ell$ particles, separated by at least $L$ vacant sites, is blocked.

The interpolation, \textit{(iii)}, is direct from the definition of the constraint $\mathbf{p}_{\ell;L}$ in \eqref{new-base} with $\ell=L$, and one can check \textit{(iv)} identically to the analogous property in Proposition \ref{prop:B-prop}. 

Property \textit{(v)} can be checked directly by decomposing $\mathbf{1}_{\{m_j^L(\eta)\geq \ell\}}=\mathbf{1}_{\{m_j^L(\eta)\geq \ell+1\}}+\mathbf{1}_{\{m_j^L(\eta)= \ell\}}$, leading to
\begin{align}
    \mathbf{p}_{\ell;L}(\eta)
    -\mathbf{p}_{\ell+1;L}(\eta)
    &=
    \frac{1}{L+1}
    \sum_{j=0}^{L}\frac{\binom{m_j^L(\eta)}{\ell}}{\binom{L}{\ell}}
    \mathbf{1}_{\{m_j^L(\eta)\geq \ell+1\}}
    \bigg(
    1-\frac{m_j^L(\eta)}{L-\ell}
    \bigg)
    \\&+
    \frac{1}{L+1}
    \sum_{j=0}^{L}\frac{\binom{m_j^L(\eta)}{l+1}}{\binom{L}{\ell+1}}
    \mathbf{1}_{\{m_j^L(\eta)=\ell\}}
    \geq 0.
\end{align}
In order to see \textit{(vi)}, one can, for example, express
\begin{align}
    \binom{m_j^L(\eta)}{\ell}
        \mathbf{1}_{
        \{
        m_j^L(\eta)
        \geq
        \ell 
        \}
        }
        =\sum_{P\in\mathcal{P}_{\ell}(W_j^{L})}
	\eta(P)
\end{align}
and note that $\int_{\Omega_N}\eta(P)\rmd\nu_\rho^N=\rho^{\ell}$ for any $P$ as above.  

\end{proof}

\addtocontents{toc}{\protect\setcounter{tocdepth}{0}}

\subsection*{Funding:}
The author gratefully acknowledges the support from FCT/Portugal through the Lisbon Mathematics PhD (Lismath), grant 47PD/BD/150345/2019; from FCT/Portugal for financial support through CAMGSD, IST-ID, projects UIDB/04459/2020 and UIDP/ 04459/2020, and through the Portugal/France agreement Programa PESSOA cotutelas (FCT/Campus France). This project was also partially supported by the ANR grant MICMOV (ANR-19-CE40-0012) of the French National Research Agency (ANR).

\addtocontents{toc}{\protect\setcounter{tocdepth}{2}}

\bibliographystyle{plain}
\bibliography{02.biblio}

\end{document}